# Quantum-assisted *hλ*-adaptive finite element method

R. H. Drebotiy[1], H. A. Shynkarenko[1]

[1]Ivan Franko National University of Lviv (Faculty of Applied Mathematics and Informatics, Department of Information Systems),

1, Universytetska St., Lviv, 79000, Ukraine

*We propose a novel finite element method scheme for singularly perturbed advection-diffusion-reaction problems, which combines certain quantum-assisted stabilization scheme with a classical h-adaptive approach to provide automatic error control and corresponding approximation refinement. Appropriate finite element a posteriori error estimates are proved. Described approach demonstrates the possibility to overcome singular perturbations by applying error-controlled smoothening to finite element approximation. Possible benefits of the proposed finite element scheme are discussed and the numerical comparison with the typical adaptive scheme is provided.*

## 1. Introduction

Advection-diffusion-reaction (ADR) model plays an important role in simulation of air and water pollution migration [17,22], optimal pollution source allocation [13] and medicine (for example, drug absorption modeling [21]). Great interest is directed to the advection-dominated problems of ADR kind. They lead us to singularly perturbed problems, whose solutions own high gradients, localized in certain parts of the domain. Approximating solutions of such problems is a complicated task [14]. Classical uniform finite element mesh refinement will not provide enough accuracy with adequate computational complexity. For such cases different kinds of mesh adaptivity are often used [3,5,12,22,23]. There are also some stabilization schemes [1,12,22], which incorporate additional penalizing terms to the initial equation to deal somehow with oscillatory nature of the obtained finite element approximations. Even with those approaches for highly dominated advection terms, it may be still needed to use huge computational resources and amount of time, needed for that computation.

Leveraging technological benefits from existing hardware is also a typical approach. In fact, most industrial finite element solvers extensively use parallelization capabilities to calculate linear system components and for solving it, providing better performance.

In addition, it may be interesting, whether it is possible to use quantum computers to obtain finite element solution. The most resource-dependent part of finite element solver is the solution of linear system. For solving linear systems on a quantum computer there are some algorithms developed. First one and most well-known is a Harrow-Hassidim-Lloyd (HHL) algorithm [15] with its improvements and variants. There are some recent developments in that direction, specially we want to note Variational Quantum Linear Solver (VQLS) quantum-classic algorithm [2], which is built to reduce the depths of the quantum circuits, needed for computations. The last part is important for the possibility of near-term usage of quantum computers for solving linear systems. It is mentioned in [8,20] for HHL algorithm, that it may be that quantum circuits for simulating unitary operators $e^{iAt}$ will have very large depth and width. It should be noted also, that all existing quantum algorithms for solving linear systems are just generating the appropriate quantum state, which is proportional to the exact solution of a system, i.e. it encodes in its amplitudes the exact solution components. Due to the quantum nature of possible measurements (i.e. wave function collapse) it is not possible in general to "read" all those components to provide the pointwise information for engineer, who uses the finite element solver. Such reconstruction is called *quantum state tomography* [19] and it can be done by making many measurements on the "equivalent" states, obtained by repeatedly using quantum linear solver. Since for finite element analysis application we operate with very large systems, such quantum tomography is impractical.

Cheaper can be to obtain some integral quantity for calculated quantum state (while it still can require several execution of the linear solver). For example, we can use so called SWAP-test [4] to calculate scalar product with some predefined state.

Detailed quantum error analysis for finite element and such a procedure with obtaining just integral quantity as a scalar product is provided in the article [18]. It shows, that HHL algorithm in general not always show exponential speedup. Also, it is important to note, that it seems, that the only efficient way to use quantum computers in the future for finite element method is calculating some integral quantity for obtained solution. At least for now, there are no known sources in literature, describing some other approaches.

In this article we propose a different strategy of possible leveraging of quantum computers for finite element method. We built certain scheme of usage of a quantum computer as an auxiliary unit for reducing high-frequency oscillations of finite element approximations for singular perturbed problems.

Mentioned strategy is based on our several previous articles [7,8,11]. Let us describe those.

In [11] we proposed new FEM stabilization scheme for 1D and 2D domains which combines regularization procedure, similar to standard Tikhonov-type regularization [16], with auxiliary Cauchy problem, corresponding to initial ADR model. Proposed scheme has one real positive regularization parameter $\lambda$, which play crucial role in obtaining accurate approximation. In [7] we proposed special algorithm for obtaining the approximate solution of this Cauchy problem, which was described in [11].

In the article [8] we proposed heuristic strategy of selecting mentioned regularization parameter also for 1D and 2D problems. Proposed approach is "quantum-assisted", since, as mentioned in [8], it can leverage quantum Harrow-Hassidim-Lloyd [15] algorithm in a combination with SWAP-test [4] to estimate regularization parameter.

In this paper we propose adaptive finite element algorithm, which combines quantum-assisted stabilization described in [8] with classical *h*-adaptivity to obtain automatic error-controlled approximation refinement algorithm. To implement that adaptivity we proved special error estimates for general variational problem (suitable for 2D and 3D) and also practical error estimate with local error indicator for 1D problem. We compare proposed algorithm with certain *h-/hp-* adaptive schemes and provide appropriate numerical experiment data.

In this article we suppose (as in our previous works), that the problem has boundary layer. Handling the case, when we have inner layers with high gradient is out of scope of current work and can be considered for future research.

We call proposed method as "$h\lambda$-adaptive", since we adapt mesh (*h*-) and the regularization parameter ($\lambda$-). Stabilization approach, which was proposed in [8] was based on linear finite elements. In this article, we also use only linear finite elements. Generalization to higher-order elements is not considered yet.

The paper is structured as follows: first we define model ADR problem; then we present and review algorithm from [11] and parameter estimation from [8]. After that we derive specialized a posteriori error estimate and describe proposed $h\lambda$-adaptive finite element scheme. In the last part we provide numerical experiment results.

## 2. Model problem

We consider the following boundary value problem (BVP) for the advection-diffusion-reaction equation:

$$\begin{cases} \text{find function } u:\bar{\Omega} \to \mathbb{R} \text{ such that:} \\ -\mu\Delta u + \vec{\beta}\cdot\nabla u + \sigma u = f \text{ in } \Omega \subset \mathbb{R}^2, \\ u = 0 \text{ on } \Gamma = \partial\Omega. \end{cases} \quad (1)$$

Here $\Omega$ is a bounded domain with a Lipschitz boundary $\Gamma = \partial\Omega$, $\mu = const > 0$ and $\sigma = const \geq 0$ are coefficients of diffusion and reaction respectively, function $f = f(x)$ and vector $\vec{\beta} = \{\beta_i(x)\}_{i=1}^2$ represent the sources and advection flow velocity respectively. We will consider noncompressible flow, i.e., $\nabla\cdot\vec{\beta} = \sum_{i=1}^2 \frac{\partial}{\partial x_i}\beta_i = 0$ in $\Omega$. By $\Delta$ we denote standard Laplace operator $\Delta = \nabla\cdot\nabla = \sum_{i=1}^2 \frac{\partial^2}{\partial x_i^2}$. Unknown function $u$ represents a substance concentration in a point of domain.

The boundary value problem (1) allows the following variational formulation:

$$\begin{cases} \text{find } u \in V := H_0^1(\Omega) \text{ such that,} \\ a(u,v) = \langle l,v \rangle \quad \forall v \in V, \end{cases} \quad (2)$$

where:

$$\begin{cases} a(u,v) = \int_\Omega (\mu\nabla u\cdot\nabla v + \vec{\beta}v\cdot\nabla u + \sigma uv)dx \quad \forall u,v \in V, \\ \langle l,v \rangle = \int_\Omega fvdx \quad \forall v \in V. \end{cases} \quad (3)$$

and $H_0^1(\Omega)$ is a standard Sobolev space with the functions having zero trace on the boundary $\Gamma$.

On $V$ we have standard Sobolev scalar product:

$$(u,v)_V = \int_\Omega (uv + \nabla u\cdot\nabla v)dx \quad (4)$$

and associated norm

$$\|u\|_V = \sqrt{(u,u)_V}. \quad (5)$$

Here and below, we assume, that components of equation (2) are quite regular and satisfy conditions of the Lax-Milgram lemma [3]. In that case, problem (2) has a unique weak solution $u \in V$, moreover, its bilinear form generates in the space of admissible functions $V$ new (energy) norm that is equivalent to the standard norm of the Sobolev space $H^1(\Omega)$.

It is well-known, that large values of the Péclet number

$$Pe := \mu^{-1}\|\vec{\beta}\|_\infty \, diam\,\Omega \quad (6)$$

where $\|\vec{\beta}\|_\infty = \left(\sum_{i=1}^2 \operatorname*{ess\,sup}_{x\in\Omega} |\beta_i(x)|^2\right)^{1/2}$, indicate, that problem (2) may be *singularly perturbed* [3,22,23]. We consider $\mu$ as a small parameter of our problem (2). Note, that in these conditions two other dimensionless quantities: Fourier and Strouhal numbers will be not large. In our consideration of the only one small parameter $\mu$, number (4) (i.e. its large values) is an indicator of singularly perturbed problem.

## 3. Stabilization scheme

To deal with large Péclet numbers, in [11] we proposed the following stabilization procedure. Let $\Gamma_0 := \{x \in \partial\Omega \mid \vec{n}(x)\cdot\vec{\beta}(x) < 0\}$, and $\vec{n}$ is a unit vector of outward normal to the boundary of domain $\Omega$. Let us introduce the following *reduced problem*:

$$\begin{cases} \text{find function } u_0 \in C^1(\Omega) \text{ such that:} \\ \vec{\beta} \cdot \nabla u_0 + \sigma u_0 = f \text{ in } \Omega, \\ u_0|_{\Gamma_0} = 0, \end{cases} \quad (7)$$

and replace the original variational problem (2) with the following *regularized problem*:

$$\begin{cases} \text{given parameter } \lambda = const \geq 0, \\ \text{given } u_0 \text{ is the solution of (7),} \\ \text{find } u^\lambda \in V \text{ such that:} \\ a(u^\lambda, v) + \lambda(u^\lambda, v)_V = \langle l, v \rangle + \lambda(u_0, v)_V \quad \forall v \in V. \end{cases} \quad (8)$$

Problem (7) is actually Cauchy problem restricted to domain $\Omega$. We assume, that there are no closed integral curves of field $\vec{\beta}$ which are entirely contained in domain $\Omega$. For solving (7) we proposed in [11] to use method of characteristics in a combination with Runge-Kutta method. Such approach can be parallelized [7] and it will not provide significant overhead to overall procedure, which we describe in this article.

**Remark 1.** Note, that in case if there are such closed integral curves, it may be, that problem (7) will be not correctly formulated. Consider, for example, equation in the form of (7) for the unknown function $w = w(x, y)$ in Cartesian coordinates:

$$\frac{-y}{\sqrt{x^2 + y^2}} w'_x + \frac{x}{\sqrt{x^2 + y^2}} w'_y = 1. \quad (9)$$

Note, that integral curves of field $\left( \dfrac{-y}{\sqrt{x^2 + y^2}}, \dfrac{x}{\sqrt{x^2 + y^2}} \right)$ are concentric circles centered at the origin. Consider, for example, integral curve $x(t) = \cos(t)$, $y(t) = \sin(t)$. Let us also consider the values of function $w(x, y)$ on that curve $z(t) = w(x(t), y(t))$. If we restrict now the equation (9) to that curve, we can clearly see, that it will transform to the equation $z'(t) = 1$. So $z(t) = t + C$, where $C = const$. Thus $w(1, 0) = w(x(0), y(0)) = z(0) = C$, but from other side: $w(1, 0) = w(x(2\pi), y(2\pi)) = z(2\pi) = 2\pi + C$. So, we have ambiguity in definition of function $w(x, y)$ by equation (9) on closed characteristic curve.

Along with problem (1) in this article we consider corresponding 1D model problem:

$$\begin{cases} \text{find function } u : \overline{\Omega} \to \mathbb{R} \text{ such that:} \\ -(\mu u')' + \beta u' + \sigma u = f \text{ in } \Omega = (0,1), \\ u(0) = u(1) = 0, \end{cases} \quad (10)$$

where $\mu = \mu(x)$, $\beta = \beta(x)$, $\sigma = \sigma(x)$, $f = f(x)$.

We provide some general error estimates, which are suited for both problems (1) and (10). For practical implementation and obtaining numerical results we derive final a posteriori error estimates only for 1D problem (10).

Corresponding variational formulation (2) for problem (10) has the following bilinear and linear forms:

$$\begin{cases} a(u, v) = \int_0^1 (\mu u' v' + \beta u' v + \sigma u v) dx, \\ \langle l, v \rangle = \int_0^1 f v dx. \end{cases} \quad (11)$$

For problem (8) we use standard finite element Galerkin formulation to obtain the approximate solution:

$$\begin{cases} \text{given parameter } \lambda = const \geq 0, \\ \text{given } u_0 \text{ is the solution of (7),} \\ \text{find } u_h^\lambda \in V_h \subset V \text{ such that:} \\ a(u_h^\lambda, v_h) + \lambda(u_h^\lambda, v_h)_V = \langle l, v_h \rangle + \lambda(u_0, v_h)_V \quad \forall v_h \in V_h. \end{cases} \quad (12)$$

where $V_h$ is a finite-dimensional subspace of $V$, consisting from piecewise-linear functions.

## 4. Stabilization parameter estimation for non-uniform mesh

In the article [8] we built certain loss function $F(\lambda)$, which indicates high-frequency parasitic oscillations in the finite element approximation. Considering problem (12) with $\lambda$ selected as such, which minimizes $F(\lambda)$, lead us to elimination of high-frequency oscillations from approximation $u_h^\lambda$. This function was built for uniform 1D meshes and arbitrary 2D meshes. Since in this article we will build adaptive scheme for 1D problem, we need proper generalization of that loss function for arbitrary 1D meshes. Suppose first, that we have *uniform* mesh with $n$ linear finite elements, which consist of points $0 = x_0 < x_1 < x_2 < ... < x_n = 1$. Let $K = [x_{i-1}, x_i]$ be the domain of $i$-th element. For finding optimal parameter $\lambda$ in [8] we proposed to use the following loss function

$$F(\lambda) = \frac{1}{|\vec{u}_h^\lambda|} \left| \sum_{k=1}^{n-2} (-1)^{k-1} \delta^2 u_h^\lambda(x_k) \right|, \tag{13}$$

where by $\delta^2 u_h^\lambda(x_k)$ we denote central difference of second order in appropriate point:

$$\delta^2 u_h^\lambda(x_k) = u_h^\lambda(x_{k+1}) - 2 u_h^\lambda(x_k) + u_h^\lambda(x_{k-1}) \tag{14}$$

and $\vec{u}_h^\lambda = (u_h^\lambda(x_1),...,u_h^\lambda(x_{n-1}))$.

Relation (13) is derived in [8] in the assumption, that the boundary layer of a problem lies near the point x=1 and thus we excluded the last element containing boundary layer from the sum (13), i.e. we have there sum only up to $k = n - 2$ excluding $\delta^2 u_h^\lambda(x_{n-1})$. Such step is precisely discussed in [8] and here we should note, that the function (13) is built to indicate the oscillations on smooth part of the solution. If we know that boundary layer is located on the left side, we can just exclude the first element instead of last. Localization of boundary layer for 2D problems can be done automatically [8].

**Remark 2.** Note, that as mentioned in [8] we can use quantum Harrow-Hassidim-Lloyd algorithm to solve the FEM linear system for certain $\lambda$ and generate quantum state $\left| \vec{u}_h^\lambda \right\rangle$ proportional to $\vec{u}_h^\lambda / |\vec{u}_h^\lambda|$. Next, we can use SWAP-test to calculate approximately the value of (13) up to a fixed multiplier constant on a quantum computer.

Expression in (13) is actually a normed scalar product between 1D discrete Laplace operator, which in some manner characterizes local curvature, calculated in the mesh points and the vector $(1, -1, 1, ..., (-1)^{n-3})$, which represents general oscillation pattern. Such "correlation" between those quantities is sensitive to the presence of similar oscillations in the approximate solution. This leads us to a function, which has unique local minimum and this minimum is reached for the optimal parameter value, for which high-frequency oscillations are excluded from the approximate solution, which is shown, in particular, by numerical experiments in [8].

In [8] we run the optimization algorithm on the $[0, \lambda_{\max}]$, where:

$$\lambda_{\max} = \frac{2 \|\beta\|_\infty \, diam\Omega}{n}. \tag{15}$$

For non-uniform meshes we need to provide suitable generalization of Laplace operator approximation (14). Here we propose to use the following modification of a quantity (9):

$$H(\lambda) = \frac{1}{|\vec{u}_h^\lambda|} \left| \sum_{k=1}^{n-2} (-1)^{k-1} [u_h^\lambda(x_{k-1}), u_h^\lambda(x_k), u_h^\lambda(x_{k+1})] \right|, \tag{16}$$

where $[u_h^\lambda(x_{k-1}), u_h^\lambda(x_k), u_h^\lambda(x_{k+1})]$ is the divided difference of order two and it is calculated in the following way:

$$[u_h^\lambda(x_{k-1}), u_h^\lambda(x_k), u_h^\lambda(x_{k+1})] := \frac{1}{x_{k+1} - x_{k-1}} \left\{ \frac{u_h^\lambda(x_{k+1}) - u_h^\lambda(x_k)}{x_{k+1} - x_k} - \frac{u_h^\lambda(x_k) - u_h^\lambda(x_{k-1})}{x_k - x_{k-1}} \right\} \quad (17)$$

Note, that quantity (16) is a direct generalization of (13) and it resembles the same generalization for 2D problems, described in [8].

## 5. A posteriori error estimates

Let us estimate the error $\| u - u_h^\lambda \|_V$. Note that we want to estimate the difference between exact solution of original problem (2) and the finite element approximation of regularized problem (8).

Let us use triangle inequality first:

$$\| u - u_h^\lambda \|_V \leq \| u - u^\lambda \|_V + \| u^\lambda - u_h^\lambda \|_V . \quad (18)$$

To estimate first term, let us subtract (2) from (8):

$$a(u^\lambda - u, v) = \lambda(u_0 - u^\lambda, v)_V . \quad (19)$$

Since we assumed, that for problem (2) the Lax-Milgram lemma conditions holds, we can write $V$-ellipticity inequality for corresponding energy norm [3]:

$$\| v \|_E^2 := a(v, v) \geq \alpha \| v \|_V^2 \quad \forall v \in V \quad (20)$$

where $\alpha = \min\{\mu, \sigma\}$ for (1). Note, that for our case, (20) also holds for 1D problem (10) with constant $\mu$ and $\sigma$.

By taking (19) with $v := u^\lambda - u$, using Cauchy-Schwarz inequality and (20) we obtain:

$$\| u - u^\lambda \|_V^2 \leq \frac{1}{\alpha} \| u - u^\lambda \|_E^2 = \frac{\lambda}{\alpha} (u_0 - u^\lambda, u^\lambda - u)_V \leq \frac{\lambda}{\alpha} \| u_0 - u^\lambda \|_V \| u - u^\lambda \|_V . \quad (21)$$

Then by dividing (21) by $\| u - u^\lambda \|_V$ and applying triangle inequality again we obtain:

$$\| u - u^\lambda \|_V \leq \frac{\lambda}{\alpha} \| u_0 - u^\lambda \|_V \leq \frac{\lambda}{\alpha} \| u_0 - u_h^\lambda \|_V + \frac{\lambda}{\alpha} \| u^\lambda - u_h^\lambda \|_V \quad (22)$$

Combining now (18) and (22) we get:

$$\| u - u_h^\lambda \|_V \leq \frac{\lambda}{\alpha} \| u_0 - u_h^\lambda \|_V + \left(1 + \frac{\lambda}{\alpha}\right) \| u^\lambda - u_h^\lambda \|_V \quad (23)$$

Consider the square of energy norm of regularized problem (8):

$$\| v \|_{E,\lambda}^2 := a(v, v) + \lambda(v, v)_V = \| v \|_E^2 + \lambda \| v \|_V^2 \geq (\alpha + \lambda) \| v \|_V^2 \quad \forall v \in V \quad (24)$$

Using (24) for second term in (23), we proved that the following estimate holds:

$$\| u - u_h^\lambda \|_V \leq \frac{\lambda}{\alpha} \| u_0 - u_h^\lambda \|_V + \frac{\sqrt{\alpha + \lambda}}{\alpha} \| u^\lambda - u_h^\lambda \|_{E,\lambda} \quad (25)$$

Taking into account elementary inequality $(a + b)^2 \leq 2(a^2 + b^2)$ we derive

**Theorem 1:** The following estimate holds:

$$\| u - u_h^\lambda \|_V^2 \leq \frac{2\lambda^2}{\alpha^2} \| u_0 - u_h^\lambda \|_V^2 + 2\frac{\alpha + \lambda}{\alpha^2} \| u^\lambda - u_h^\lambda \|_{E,\lambda}^2 \quad (26)$$

Note, that (26) holds for the problems of kind (1) and (8) in any number of dimensions.

All quantities under the norm $\| u_0 - u_h^\lambda \|_V$ are known, so this norm can be approximately computed by quadrature.

The second term norm $\| u^\lambda - u_h^\lambda \|_{E,\lambda}$ measures the finite element error for regularized problem in its natural energy norm, which makes it suitable for applying some known error estimators.

Let us consider the 1D boundary value problem, which corresponds to the regularized problem (8). Consider for small simplification the case $\mu = const$ (it is quite practical in actual applications). Using integration by parts, it is not hard to obtain the following problem:

$$\begin{cases} \text{find function } u^\lambda : \overline{\Omega} \to \mathbb{R} \text{ such that:} \\ -(\mu+\lambda)(u^\lambda)'' + \beta \cdot (u^\lambda)' + (\sigma+\lambda)u^\lambda = f + \lambda u_0 - \lambda u_0'' \text{ in } \Omega = (0,1), \\ u^\lambda(0) = u^\lambda(1) = 0, \end{cases} \tag{27}$$

Taking into account, that we use linear elements, i.e. second derivative of approximation is almost everywhere zero, we can define the residual:

$$R[u_h^\lambda] := f + \lambda u_0 - \lambda u_0'' - \beta(u_h^\lambda)' - (\sigma+\lambda)u_h^\lambda \tag{28}$$

Let us define quadratic bubble-function $\omega_K(x) := (x_i - x)(x - x_{i-1})$.

Let us now use the explicit residual error estimator described in [9] (which is a direct generalization of estimator from [5] to ADR problems) for our case of linear elements:

$$\| u^\lambda - u_h^\lambda \|_{E,\lambda}^2 \leq \sqrt{2} \left( \sum_K \| \sqrt{\omega_K} R[u_h^\lambda] \|_{L^2(K)}^2 \right)^{\frac{1}{2}} \| u^\lambda - u_h^\lambda \|_V \leq$$
$$\leq \sqrt{\frac{2}{\alpha+\lambda}} \left( \sum_K \| \sqrt{\omega_K} R[u_h^\lambda] \|_{L^2(K)}^2 \right)^{\frac{1}{2}} \| u^\lambda - u_h^\lambda \|_{E,\lambda} \tag{29}$$

or finally:

$$\| u^\lambda - u_h^\lambda \|_{E,\lambda}^2 \leq \sum_K \left\{ \sqrt{\frac{2}{\alpha+\lambda}} \| \sqrt{\omega_K} R[u_h^\lambda] \|_{L^2(K)} \right\}^2 \tag{30}$$

Let us introduce per-element Sobolev norm by the expression:

$$\| v \|_K^2 := \int_K \left( (v')^2 + v^2 \right) dx \tag{31}$$

Note, that:

$$\| v \|_V^2 = \sum_K \| v \|_K^2 \tag{32}$$

Now, using (30), we can rewrite (26) in the following form:

$$\| u - u_h^\lambda \|_V^2 \leq \sum_K \left\{ \frac{2\lambda^2}{\alpha^2} \| u_0 - u_h^\lambda \|_K^2 + \frac{4}{\alpha^2} \| \sqrt{\omega_K} R[u_h^\lambda] \|_{L^2(K)}^2 \right\} = \sum_K \eta_K^2 \tag{33}$$

where

$$\eta_K := \frac{1}{\alpha} \left\{ 2\lambda^2 \| u_0 - u_h^\lambda \|_K^2 + 4 \| \sqrt{\omega_K} R[u_h^\lambda] \|_{L^2(K)}^2 \right\}^{\frac{1}{2}} \tag{34}$$

Quantity $\eta_K$ will be used in the next chapter as a local error indicator for the adaptivity step.

## 6. *hλ*-adaptivity

Let us bring all obtained pieces together and describe the quantum-assisted *hλ*-adaptive finite element method. Consider the auxiliary procedure:

**Algorithm 1.** *Procedure* **QuantumStabilizedSolution**($\mathfrak{I}_h$, $\lambda_{max}$)

**Require:** mesh $\mathfrak{I}_h$; max parameter bound $\lambda_{max}$;

1. Calculate stabilized approximation $u_h^\lambda$ and corresponding optimal parameter $\lambda^*$ by **Algorithm 2 from [8]** with usage of loss function $H(\lambda)$ defined by (16) and using max parameter bound $\lambda_{max}$.
2. **return pair(** $u_h^\lambda$, $\lambda^*$ **)**

**Remark 3.**
Since we use relatively coarse estimate (33), for large advection term we will have very small values of obtained concentration. In general, to see meaningful results, we need to set large enough source term, but in that case the residual terms in (33) will be very large, leading to very large total value of error estimate. To overcome this, we propose to use for total error quantity relative-to-source error estimate, by dividing estimate (33) by some norm of source term.

**Remark 4.**
As in [8], we do not expect the approximation of solution in the zone of boundary layer, that's why for total error control we will use total error estimate excluding boundary layer. Given mesh $\mathfrak{I}_h$ and boundary layer separation node $x_{layer}$ let us define two sets of finite elements:

$$Left := elements\left(\{x_i \in \mathfrak{I}_h \mid x_i \leq x_{layer}\}\right) \quad (35)$$

$$Right := elements\left(\{x_i \in \mathfrak{I}_h \mid x_i \geq x_{layer}\}\right) \quad (36)$$

By taking into account **Remark 3**, we define:

$$Error = \sqrt{\frac{\sum_{K \in Left} \eta_K^2}{\sum_{K \in Left} \| f \|_{L^2(K)}^2}} \times 100\% \quad (37)$$

**Algorithm 2.** *Quantum-assisted hλ-adaptive FEM* (*N, Tol,* $\theta$)

**Require:**

    initial mesh element count *N*;

    acceptable relative-to-source error level *Tol*;

    refinement threshold $\theta$.

1. **Initialization:** generate initial uniform *N*-element mesh $\mathfrak{I}_h$,
   assign the pre-last node to $x_{layer}$,
   calculate $\lambda_{max}$ using formula (15)
2. $u_h^\lambda$, $\lambda^*$ := **QuantumStabilizedSolution**($\mathfrak{I}_h$, $\lambda_{max}$)
3. Calculate Error using formula (37):
4. **If** Error < Tol **then return** $u_h^\lambda$
5. Calculate error indicators $\eta_K$ using formula (34) for $K \in Left \cup Right$
6. Refine current mesh $\mathfrak{I}_h$:
   Bisect all elements $K \in Left$ for which $\eta_K > \theta \max_{K \in Left} \eta_K$
   Bisect all elements $K \in Right$ for which $\eta_K > \theta \max_{K \in Right} \eta_K$

7. **goto** step 2.

**Remark 5.**

One way of possible optimization of the above algorithm is passing previous $\lambda^*$ to the stabilization algorithm as $\lambda_{max}$ on Step 2. Such $\lambda_{max}$ reassignment provides some optimization of the number of quantum solver calls, since theoretically we can expect that optimal regularization parameter value will decrease on further mesh refinement, since the approximation capability of the FE spaces is increased. In practice, when the initial mesh does not provide enough precision, it may be not the case and subsequent optimal lambdas can be not always monotonically decreasing.

**Remark 6.**

Note the way of doing refinements in the algorithm. We do refinements separately for zones outside and inside the boundary layer. This is due to the need of approximate the solution precisely outside the boundary layer and avoid too extensive boundary layer refinement.

## 7. Numerical experiment

Let us consider 1D singularly perturbed ADR BVP with the following data:

$$\mu = 1, \beta = 10^4, \sigma = 10^2, f = 10^4 \cos(4.5\pi x) \tag{38}$$

For the algorithm let us set $\theta = 2/3, Tol = 0.2\%, N = 15$.

Algorithm stopped after 5 iterations. On Fig. 1 obtained approximations are depicted. Since real quantum computer with needed qubit number and error-tolerant routines is not available for us, we executed the part for quantum computer (calculating $H(\lambda)$ by solving linear system to obtain stabilized approximation) on the classic computer. Number of needed calls to quantum computer circuit is shown in the Table 1. By the dashed line, we show exact solution approximation obtained on uniform mesh with 5000 finite elements and on the horizontal axis we show distribution of finite elements.

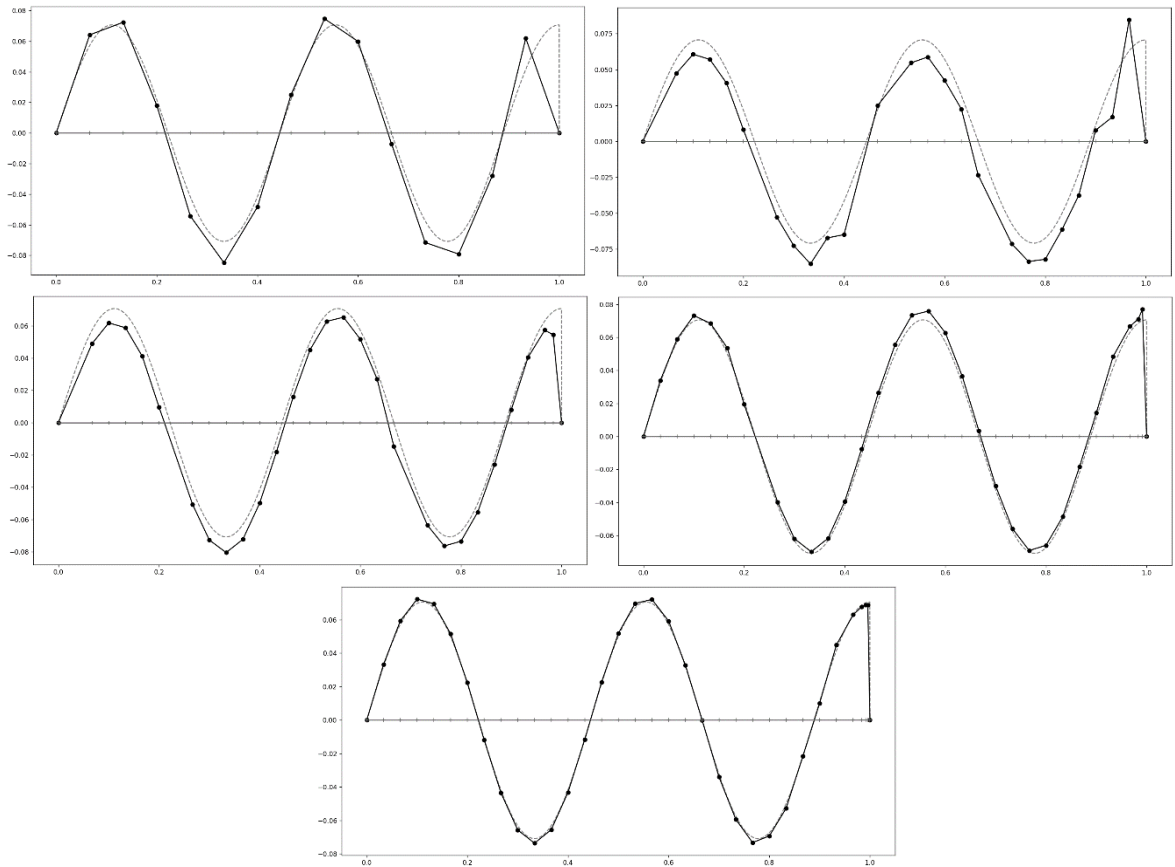

**Fig. 1.** *Plots of subsequent approximations of hλ-adaptive scheme, ordered from left to right, top to bottom.*

Let us consider absolute error estimator by the formula $\eta = \sqrt{\sum_{K \in Left} \eta_K^2}$. Suppose, that we have two subsequent iterations of the algorithm with the corresponding numbers of d.o.f. $N_1, N_2$ and computed absolute error estimates $\eta_1, \eta_2$. We can estimate order of convergence in the following way:

$$\{Order\ of\ convergence\} = -\frac{\ln \eta_2 - \ln \eta_1}{\ln N_2 - \ln N_1}. \tag{39}$$

In the next table we present computed convergence results.

**Table 1.** *Convergence results for hλ-adaptive scheme.*

| Iteration # | D.o.f. count | Accumulated d.o.f. count | Error, % | Order of convergence | Number of needed quantum solver calls (HHL+SWAP) |
|---|---|---|---|---|---|
| 1 | 14 | 14 | 2.32 | - | 31 |
| 2 | 24 | 38 | 1.46 | 0.85 | 35 |
| 3 | 27 | 65 | 0.90 | 4.17 | 33 |
| 4 | 30 | 95 | 0.32 | 9.68 | 37 |
| 5 | 32 | 127 | 0.16 | 10.78 | 19 |

**Remark 7.** For comparison we run classical *h*-adaptivity with the same error indicator ((34) with λ=0). Comparable results were obtained in 8 iterations with total 541 d.o.f. used on all iterations. On the Fig. 2 we show the corresponding plots of subsequent approximations. In table 2 we present convergence history.

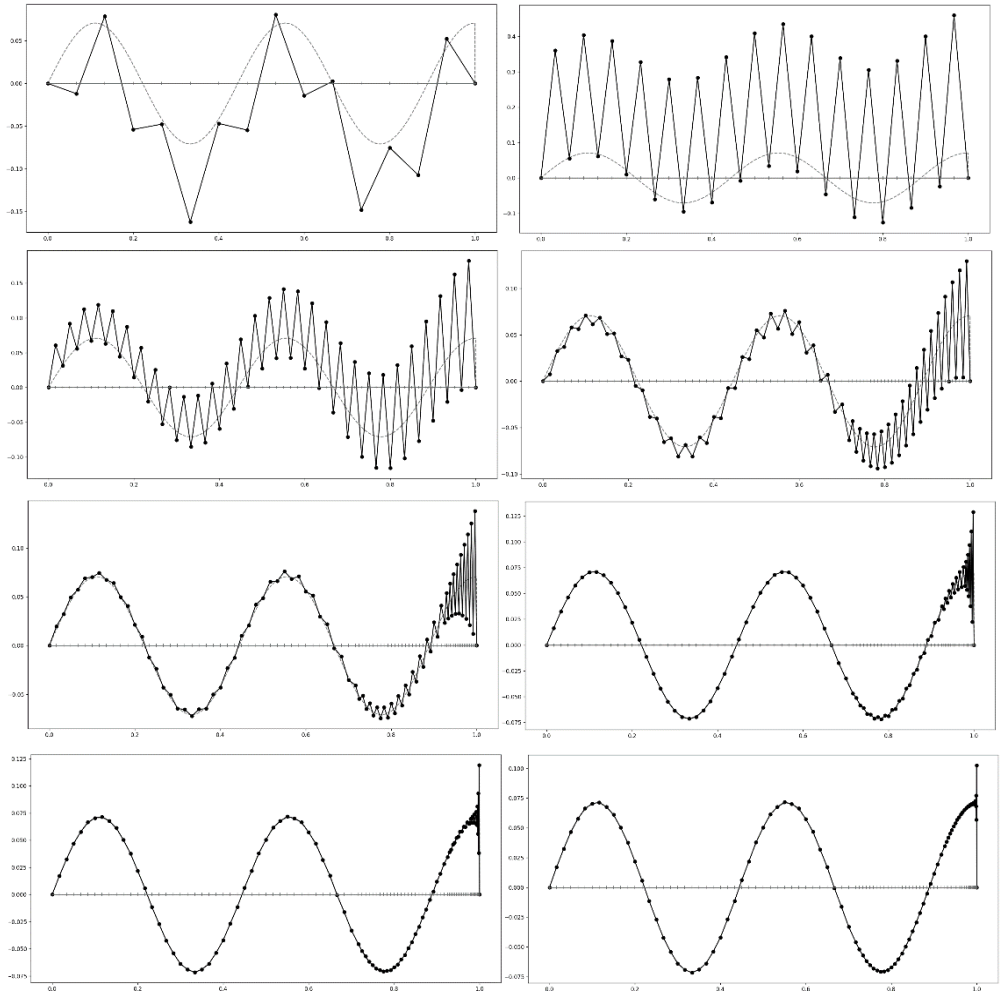

**Fig. 2.** *Plots of subsequent approximations of h-adaptive scheme, ordered from left to right, top to bottom.*

**Table 2.** *Convergence results for h-adaptive scheme.*

| Iteration # | d.o.f. count | Accumulated d.o.f. count | Error, % | Order of convergence |
|---|---|---|---|---|
| 1 | 14 | 14 | 9.33 | - |
| 2 | 29 | 43 | 44.25 | -2.14 |
| 3 | 59 | 102 | 11.25 | 1.93 |
| 4 | 76 | 178 | 3.40 | 4.72 |
| 5 | 85 | 263 | 1.06 | 10.39 |
| 6 | 90 | 353 | 0.21 | 28.57 |
| 7 | 93 | 446 | 0.09 | 24.01 |
| 8 | 95 | 541 | 0.07 | 15.63 |

In addition, we run *hp*-adaptive scheme described in [9,10] with the error indicator from [6] based on the local fundamental solutions. Maximum polynomial degree used – 5. In Fig. 3 we depict all subsequent approximations and in the table 3 we show convergence history.

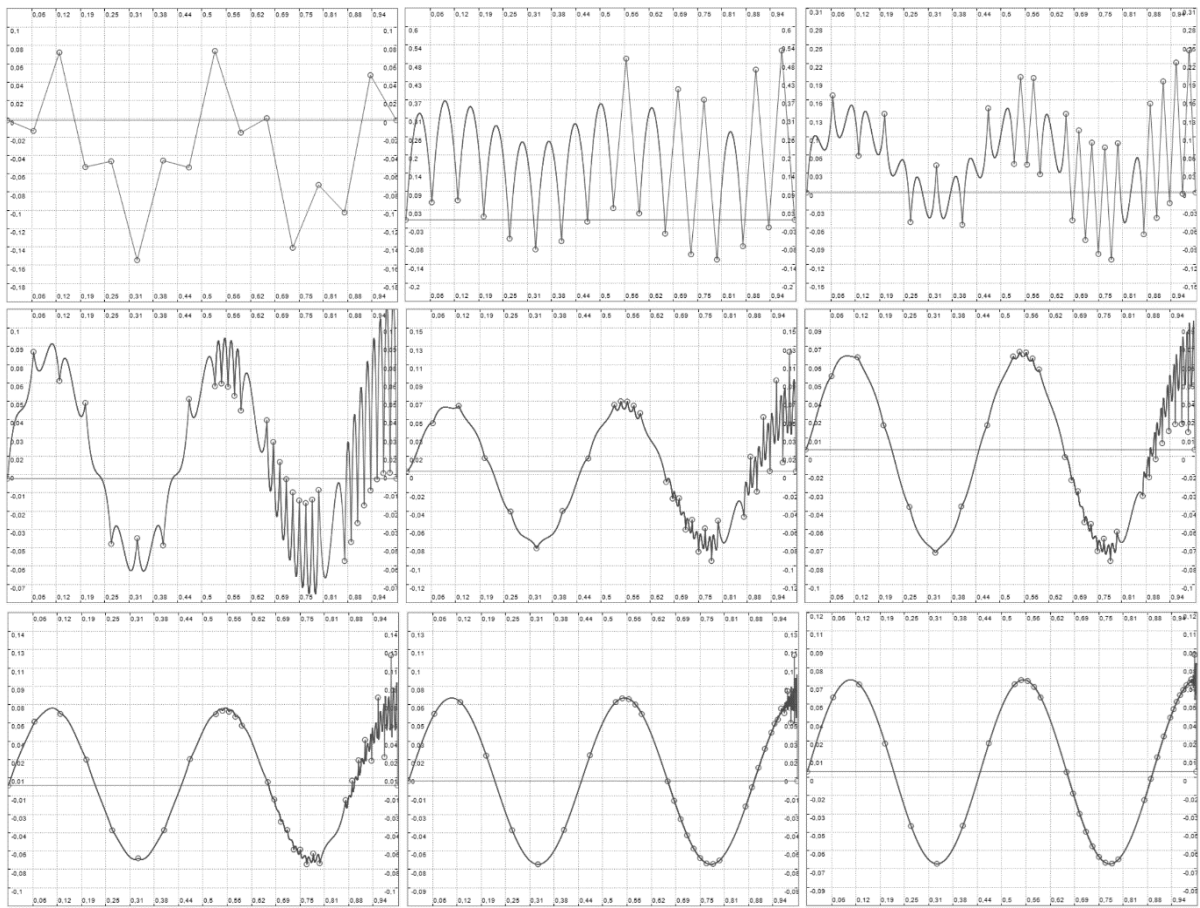

**Fig. 3.** *Plots of subsequent approximations of hp-adaptive scheme, ordered from left to right, top to bottom.*

**Table 3.** *Convergence results for hp-adaptive scheme.*

| Iteration # | d.o.f. count | Accumulated d.o.f. count | Error, % | Order of convergence |
|---|---|---|---|---|
| 1 | 14 | 14 | 61.74 | - |
| 2 | 29 | 43 | 127.45 | -0.99 |
| 3 | 49 | 92 | 183.95 | -0.69 |

| 4 | 69  | 161 | 629.67 | -3.59 |
|---|-----|-----|--------|-------|
| 5 | 84  | 245 | 117.03 | 8.55  |
| 6 | 92  | 337 | 84.72  | 3.55  |
| 7 | 97  | 434 | 90.53  | -1.25 |
| 8 | 117 | 551 | 23.64  | 7.16  |
| 9 | 127 | 678 | 14.03  | 6.35  |

First, note, that we can see different scales for the error in tables 1-3 due to different estimators used (for *h*- vs *hp*-) or due to the stabilization used which adds additional quantity to the estimator (for *hλ*- vs *h*-). Note also, that total d.o.f. and iteration count is less for *hλ*- scheme as comparing to *h*- and *hp*-.

We see clearly meaningful shape of the approximation yet from first iteration for the *hλ*-adaptive scheme on Fig. 1. This is not true for classical *h*-adaptivity and *hp*-adaptivity also, since first iterations are generating very oscillatory approximations, making them not usable for qualitative analysis of the solution by the researcher. With such high gradient in the boundary layer as we have in considered model problem, oscillations spread over entire solution domain. To approximate the solution in *h*- and *hp*- schemes precisely, we need to approximate it in the boundary layer, so in that case using mesh refinements (*h*-) is even more effective than adding higher polynomial degree (*p*-), for which we clearly observe even non-monotonic convergence.

**Remark 8.** We should note, that algorithm of parameter estimation will not work properly, if the initial number of elements is too low. In that case approximation of Laplace operator by the divided differences is not enough accurate in the sense, that for smoothed approximation we will still have large enough difference between two values, computed for adjacent mesh nodes, breaking the expected shape of the function $H(\lambda)$ (or $F(\lambda)$). How to set initial mesh element number is an open question. At least, it is reasonable to take it comparable to the one, on which similar but not singularly perturbed problem will be solved adequately on the uniform mesh.

**Remark 9.** According to [8] the number of quantum solver calls will be bounded by some constant, which depends on a precision for parameter estimation (which in presented example was set to 1) and a value of $\lambda_{max}$.

## 8. Conclusions

In this paper we proposed quantum-assisted *hλ*-adaptive finite element method scheme for singularly perturbed advection-diffusion-reaction boundary value problems. We provided heuristic criteria for regularization procedure for non-uniform meshes for 1D domains, based in similar generalization for 2D domains. Also, we derived certain a posteriori error estimates for general case (suited for 1D, 2D, 3D) and special estimator with local error indicators for 1D. Numerical experiment data is included and algorithm was compared to classical certain *h*-/*hp*-adaptive schemes with. Despite that we implemented practical algorithm only for 1D domains, proposed method can be generalized to higher dimensions and, probably, to different types of boundary value problems, but this is out of scope of current paper.